\documentclass[12pt]{article}

\usepackage{amsfonts}

\begin{document}

\begin{center}
 {\large \bf Inverse scattering problem for Sturm--Liouville operator on non-compact A--graph. Uniqueness result.
 } \\[0.2cm]
 {\bf Mikhail Ignatyev} \\[0.2cm]
\end{center}

\thispagestyle{empty}

{\bf Abstract.} In a finite--dimensional Euclidian space we consider a connected metric graph with the following property: 
each two cycles can have at most one common point. Such graphs are called A-graphs. 
On noncompact A-graph we consider a scattering problem for Sturm--Liouville
differential operator with
standard matching conditions in the internal vertices. Transport, spectral and scattering problems
for differential operators on graphs appear frequently in mathematics,
natural sciences and engineering. In particular, direct and inverse problems
for such operators are used to construct and study models in mechanics, nano-electronics,
quantum computing and waveguides. The most complete results on (both direct and
inverse) spectral problems were achieved in the case of
Sturm--Liouville operators on compact graphs, in the noncompact case there are no similar general results.
In this paper, we establish some properties of the spectral characteristics and investigate the inverse problem of recovering the operator from the scattering data.
A uniqueness
theorem for such inverse problem is proved.\\
[0.1cm]



{\bf 1. Introduction.}
Let $\Gamma$ be a metric graph in a finite--dimensional Euclidian space
$\mathcal{H}$ with the set of vertices
$V(\Gamma)$ and the set of edges $\mathcal{E}(\Gamma)\cup\mathcal{R}(\Gamma)$, where $\mathcal{E}(\Gamma)$ is the set of compact edges and $\mathcal{R}(\Gamma)$ is the
set of rays. We assume that all edges are the smooth curves in $\mathcal{H}$ which can intersect only in the vertices. We parameterize all the edges with the natural parameters and for 2 any points $x, x'$ of the same edge we denote $|x-x'|$ the distance between these points along the edge (i.e., the corresponding arc length).
Let $y(\cdot)$ be some function on $\Gamma$. For $x\in\mbox{int}\ r$,
$r\in\mathcal{E}\cup\mathcal{R}$ we define $y'(x)$ as the derivative
with respect to the local natural parameter (i.e. arc length) on the edge $r$. Thus, we can determine the Laplacian operator $y''(x)$ for $x\in\mbox{int}\ r$. Then, for $v\in V$ and $r\in \mathcal{E}(\Gamma)\cup\mathcal{R}(\Gamma)$ such that $r$ is incident to $v$ we define
$\partial_r y(v)$ as the derivative in  direction to
the interior of the edge $r$. We denote by $I(v)$ the set of all edges which are
incident to $v$. The following condition in the internal
vertex $v$ is called {\it standard matching condition} and is
denoted $MC(v)$:
$$
\sum\limits_{r\in I(v)}\partial_r y(v)=0. \eqno(1.1)
$$
Suppose that the set $\partial\Gamma$ of boundary vertices is divided into 2 parts: $\partial\Gamma=\partial_K\Gamma\cup\partial_D\Gamma$. We call the vertices from $\partial_D\Gamma$ as {\it $D$-type vertices} and the vertices from $\partial_K\Gamma$ as {\it $K$-type vertices}.
For $K$-type vertices we assume the $MC(v)$ condition in the form (1.1) (that obviously becomes the Neumann condition), for $v\in\partial_D\Gamma$ we use the Dirichlet condition:
$$
y(v)=0 \eqno(1.2)
$$
in the capacity of $MC(v)$.

Now let $q(x)$ be a real-valued integrable
function on $\Gamma$ satisfying the following condition:
$$\int\limits_r(1+|x|)|q(x)|\;d|x|<\infty \eqno(1.3)$$
for all $r\in\mathcal{R}$, where $|x|$ is a natural parameter on $r$ measured from the initial point of the ray.
We consider
the differential expression
$$
\ell y:=-y''+q(x)y \eqno(1.4)
$$
and the Sturm--Liouville operator $L=L(\Gamma,q)$ in $L_2(\Gamma)\cap C(\Gamma)$ which
is generated by the expression (1.4) and the matching conditions
$MC(v), v\in V$. More exactly, we assume $y\in dom L$ iff $y\in L_2(\Gamma)\cap C(\Gamma)$, $y$ belongs to $W_2^2(r)$ for each $r\in\mathcal{E}\cup\mathcal{R}$ and satisfies the matching conditions $MC(v)$ for all $v\in V$.

Transport, spectral and scattering problems
for differential operators on graphs appear frequently in mathematics,
natural sciences and engineering {\cite{FaP}, \cite{Ex}, \cite{KoS},
\cite{Nic}, \cite{PoB}, \cite{PoP}}. In particular, direct and inverse problems
for such operators are used to construct and study models in mechanics, nano-electronics,
quantum computing and waveguides {\cite{Kuch}, \cite{AKu}}.

During the last years such problems were in the focus of intensive
investigations. The most complete results on (both direct and
inverse) spectral problems were achieved in the case of
Sturm--Liouville operators on compact graphs {\cite{BrW},
\cite{Bel1}, \cite{Yur1}, \cite{Yur2}, \cite{Yur3}, \cite{Yur4},
\cite{Yur5}}, where certain systems of spectra or Weyl functions
where  shown to be appropriate input data for the inverse problems
and where also some constructive procedures for solving these
problems were developed.

In the noncompact case there are no similar general results. The
presence of several noncompact edges ({\it rays}) and compact
edges simultaneously leads to some qualitative difficulties in the
investigation of the spectral problems due to the non-classical
behavior of the main objects, such as Weyl--type solutions and
reflection coefficients for the rays. For the first time the
scattering problem on noncompact graphs was considered
systematically in {\cite{Ger}}, where some useful observations
were made, but complete results have been obtained only for the
special case of star--type graphs. In {\cite{FrY1}} an inverse
spectral problem on noncompact graphs with one ray has been
investigated using Weyl functions. In {\cite{Troo}} the authors
solved a particular inverse scattering problem of recovering an
operator on the ray of the simplest noncompact graph consisting of
one cycle and one ray. Some results for graphs consisting of one cycle and several rays were obtained in {\cite{FrI}}, {\cite{ISI}}. One should also mention the works
{\cite{KuS}, \cite{KuB}}, where some non-uniqueness results where
obtained for inverse scattering problems on general noncompact
graphs.

In this paper we study the Sturm--Liouville operators on connected noncompact graphs with the following property:
each two cycles can have at most one common point (here and everywhere in this paper cycle is a chain of different edges that forms a closed curve). Such graphs are called A-graphs. Main result of the paper is the uniqueness theorem for the inverse scattering problem (see Theorem 8.1).

We complete this section with some technical remarks. Let $G$ be some A-graph. The edge $r$ is called {\it simple edge} if it is not part of any cycle. We agree to call simple edges and cycles as {\it $a$-edges}, i.e. $a$-edge is either cycle or simple edge.  If some simple edge is incident to the boundary vertex then the edge is called {\it boundary edge}. For definiteness we assume that $G$ has at least one boundary vertex (this assumption is not necessary for the assertion of Lemma 2.7 and Corollary 2.4), let us take one of them as a root. We denote it as $v^0$ and the corresponding boundary edge as $r^0$.
Let $\mathfrak{a}$ be some $a$-edge. The minimal number $\omega_\mathfrak{a}$ of $a$-edges between the rooted edge and $\mathfrak{a}$ (including $\mathfrak{a}$) is called the {\it order} of $\mathfrak{a}$. The order of rooted edge is equal to zero. Let $\mathcal{C}$ be the set of all cycles and $\mathcal{A}$ be the set of all $a$-edges. The number $\omega:=\max\limits_{\mathfrak{a}\in\mathcal{A}}\omega_\mathfrak{a}$ is called the order of graph $G$. We denote by $\mathcal{A}^{(\mu)}$, the set of $a$-edges of order $\mu$.

For given $\mathfrak{a}\in\mathcal{A}$ we define the graph $G^+(\mathfrak{a})$ as a union of all the edges $r$ with the property: any path containing $r$ and the rooted edge $r^0$ necessarily contains some edge from $\mathfrak{a}$. For given cycle $\mathfrak{c}\in\mathcal{C}$ we define the graph $G_\mathfrak{c}$ as follows. Let $\mathfrak{c}$ consists of the edges (subsequently) $r_1, r_2,\ldots, r_p$ that connecting $v_0$ with $v_1$, $v_1$ with $v_2$, ..., $v_{p-1}$ with $v_0$, where $v_0=:u_\mathfrak{c}$ is the vertex of $\mathfrak{c}$ nearest to the root. Then $G_\mathfrak{c}$ is a graph obtained from $G$ by replacing the edge $r_p$ connecting $v_{p-1}$ and $v_0$ with the edge $r'_p$ of the same length connecting $v_{p-1}$ and some additional vertex $v_\mathfrak{c}$. We identify the points of $r'_p$ with the points of $r_p$ and set $q\left|_{r'_p}\right.=q\left|_{r_p}\right.$, i.e. we set $q(x', G_\mathfrak{c}):=q(x,G)$ for $x'\in r'_p$, $x\in r_p$ such that $|x'-v_{p-1}|=|x-v_{p-1}|$. We call $\mathfrak{c}$ a {\it boundary cycle} if $G^+(\mathfrak{c})=\mathfrak{c}$.

\medskip
{\bf 2. Eigenvalue problem. Weyl functions and characteristic functions.} Let $\Gamma$ be an arbitrary noncompact (not necessarily connected) graph. It is well-known that an eigenvalue problem for $L=L(\Gamma, q)$ can be reduced to some linear algebraic system by using the local FSS on each edge.

Suppose that we choose and fix some (arbitrary) ordering $\prec$ on $V$, $\mathcal{E}$ and $\mathcal{R}$ (we agree to use the same symbol for all these three orderings). For the edge $r\in\mathcal{E}$ connecting 2 vertices $u$ and $v$, where $u\prec v$ we agree to consider $u$ as an initial vertex and $v$ as a terminal vertex for $r$. Also we shall use the notation $r=[v,\infty)$ for the ray $r$ emanating from $v$. Let $y(\cdot)$ be the eigenfunction corresponding to the eigenvalue $\lambda\in\mathbf{C}\setminus [0,+\infty)$. We write $y(\cdot)$ in the following form:
$$
y(x)=\beta^1_r C_r(x,\lambda)+\beta^2_r S_r(x,\lambda), x\in r\in
\mathcal{E}. \eqno(2.1)
$$
$$
y(x)=\gamma_r e_r(x,\rho), x\in r\in \mathcal{R}, \eqno(2.2)
$$
Here $C_r(x,\lambda), S_r(x,\lambda)$ are the cosine- and sin- type
solutions for the equation $\ell y=\lambda y$ on the edge $r$,
i.e., $C_r(x,\lambda), S_r(x,\lambda)$ are satisfying the initial conditions:
$$
C_r(u,\lambda)=\partial_r S_r(u,\lambda)=1, \quad S_r(u,\lambda)=\partial_r C_r(u,\lambda)=0,
$$
where $u,v\in V$, $u,v\in r$, $u\prec v$. Then, $e_r(x,\rho)$ is the Jost
solution for the equation $\ell y=\lambda y$ on the ray $r$, $\lambda=\rho^2$, $\rho\in\Omega_+:=\{\rho:\ \mbox{Im}\rho>0\}$.

In view of (2.1) and (2.2), the
matching conditions $MC(v), v\in V$ together with the condition
$y(\cdot)\in C(\Gamma)$ reduce to a system of linear algebraic
equations with respect to the values $\{\beta^1_r,
\beta^2_r\}_{r\in\mathcal{E}}$, $\{\gamma_r\}_{r\in\mathcal{R}}$ and $\alpha_u:=y(u), u\in V$. More precisely, we assign each compact edge $r$, $u,v\in r$, $u,v\in V$ with the following pair of equations:
$$
\beta^1_r C_r(u,\lambda)+\beta^2_r S_r(u,\lambda)-\alpha_u=0, \ \beta^1_r C_r(v,\lambda)+\beta^2_r S_r(v,\lambda)-\alpha_v=0, \eqno(2.3)
$$
each ray $r=[v,\infty)$ with the equation:
$$
\gamma_r e_r(v,\rho)-\alpha_v=0, \eqno(2.4),
$$
each internal vertex $v$ and each vertex $v\in\partial_K \Gamma$ with the equation:
$$
\sum\limits_{r\in I(v)\cap \mathcal{E}} \left(\beta^1_r \partial_r C_r(v,\lambda)+\beta^2_r \partial_r S_r(v,\lambda)\right)+
\sum\limits_{r\in I(v)\cap \mathcal{R}} \gamma_r \partial_r e_r(v,\rho)=0 \eqno(2.5)
$$
and each vertex $v\in\partial_D \Gamma$ with the equation:
$$
\alpha_v=0. \eqno(2.6)
$$
We set the ordering of equations in group (2.3), (2.4) according to the ordering of edges and ordering of equations (2.5), (2.6) according to the ordering of vertices. Now we define the characteristic function $\Delta(\lambda, \Gamma, q)$ as a determinant of the system (2.3) -- (2.6) with the ordering of equations described above. One can notice that, for given graph $\Gamma$ and potential $q(\cdot)$, the characteristic function is uniquely determined by the ordering of the edges and vertices. Below, if the potential $q$ is the same for all the terms in the relation, we shall often omit $q$ in argument's list of $\Delta$ and write $\Delta(\lambda,\Gamma)$ instead of $\Delta(\lambda, \Gamma,q)$.

Let $v$ be an arbitrary internal vertex or $K$-type boundary vertex. Denote $E(\Gamma, v)$ the graph
$\Gamma_0$ which is constructed in the following way: \\[5mm]
1) replace the vertex $v$ with the set of vertices $\{v'_r\}_{r\in I(v)}$;\\
2) replace each $r\in I(v)$ connecting $v$ and some vertex $u$ with the edge
of the
same length connecting $u$ and $v'_r$ , replace $r=[v,\infty)\in\mathcal{R}(\Gamma)$ with the ray $[v'_r,\infty)$; \\
3) add no other vertices or edges. All the $v'_r, r\in I(v)$
become boundary vertices of $\Gamma_0$, we assume that all of them are $D$-type vertices.\\
We assume that the additional vertices $v'_r, r\in I(v)$ are ordered according to the ordering of corresponding edges $r\in I(v)$.
Clearly we can assume that we identify
$\mathcal{E}(\Gamma_0)$ with $\mathcal{E}(\Gamma)$ and
$\mathcal{R}(\Gamma_0)$ with $\mathcal{R}(\Gamma)$ and consider the
Sturm--Liouville operator $L(\Gamma_0,q)$ with the same
potential $q(\cdot)$.

\medskip
Let us take an arbitrary vertex $v\in V(\Gamma)$. We call a function
$\Phi_v(x,\lambda,\Gamma)$, defined at least for $x\in\Gamma$, $\lambda\in
\mathbf{C}\setminus\mathbf{R}$ {\it the Weyl solution associated
with $v$} iff:\\[5mm]
1) it is continuous on $\Gamma$ (with respect to $x$) and satisfying $MC(u)$ for all $u\in V(\Gamma)\setminus\{v\}$; \\
2) it solves the differential equation $\ell \Phi_v=\lambda \Phi_v$, $x\in
\mbox{int}\
r$, $r\in \mathcal{E}(\Gamma)\cup\mathcal{R}(\Gamma)$; \\
3) $\Phi_v(\cdot,\lambda,\Gamma)\in L_2(\Gamma)$;\\
4) $\Phi_v(v,\lambda,\Gamma)=1$.\\[5mm]
The value
$$M_v(\lambda,\Gamma):=\sum\limits_{r\in I(v)}\partial_r\Phi_v(v,\lambda,\Gamma)$$
is called {\it the Weyl function associated with $v$}.\\

We recall some facts concerning the properties of Weyl functions {\cite{FrI}}.

\medskip
{\bf Lemma 2.1.} {\it $M_v(\lambda,\Gamma)$ is a Nevanlinna function.}

\medskip
{\bf Lemma 2.2.} {\it For the Weyl function $M_v(\lambda,\Gamma)$ associated
with the internal vertex or $K$-type boundary vertex $v$ the following representation holds:
$$ M_v(\lambda, \Gamma)=\frac{\Delta(\lambda)}{\Delta_0(\lambda)}, $$
where $\Delta(\lambda)=\Delta\left(\lambda, \Gamma\right)$,
$\Delta_0(\lambda)=\Delta\left(\lambda, \Gamma_0\right)$ and
$\Gamma_0=E(\Gamma,v)$.}

\medskip
{\bf Lemma 2.3.} {\it Suppose that the graph $\Gamma$ is represented as $\Gamma=\Gamma_1\cup\Gamma_2$, where the graphs $\Gamma_1$, $\Gamma_2$ are such that $\Gamma_1\cap\Gamma_2=v$. Let the ordering $\prec$ on $\Gamma_1$, $\Gamma_2$ be inherited from the ordering on $\Gamma$. If, further, $v$ is boundary vertex for $\Gamma_j$ we assume it to be a $K$-type vertex.\\ Then the following representation holds:
$$
\Delta(\lambda,\Gamma)=\Delta(\lambda,\Gamma_1)\Delta(\lambda,\Gamma'_2)+\Delta(\lambda,\Gamma'_1)\Delta(\lambda,\Gamma_2),
$$
where $\Gamma'_j=E\left(\Gamma_j,v\right)$.}

\medskip
{\bf Proof.} Using the representation from Lemma 2.2 we write:
$$
M_v(\lambda,\Gamma)=\frac{\Delta(\lambda,\Gamma)}{\Delta(\lambda,\Gamma')},
$$
where $\Gamma'=E\left(\Gamma,v\right)$. Since $\Gamma'=\Gamma'_1\cup\Gamma'_2$, $\Gamma'_1\cap\Gamma'_2=\emptyset$ we have $\Delta(\lambda,\Gamma')=\Delta(\lambda,\Gamma'_1)\Delta(\lambda,\Gamma'_2)$. Further, it's clear that $\Phi_v(x,\lambda,\Gamma)=\Phi_v(x,\lambda,\Gamma_j)$ for $x\in\Gamma_j$. Thus, one can easily obtain that $M_v(\lambda,\Gamma)=M_v(\lambda,\Gamma_1)+M_v(\lambda,\Gamma_2)$. Substituting here the representations from Lemma 2.2 we obtain the required relation.$\hfil\Box$

\medskip
{\bf Corollary 2.1.} {\it Let $\Gamma=\bigcup\limits_{j=1}^p\Gamma_j$, where $\Gamma_j\cap\Gamma_k=v$ for any $j\neq k$ and $v\in\partial_K\Gamma_j$ if $v\in\partial\Gamma_j$ . Then
$$
\Delta(\lambda,\Gamma)=\sum\limits_{k=1}^p\Delta(\lambda,\Gamma_k)\prod\limits_{j\neq k}\Delta(\lambda,\Gamma'_j),
$$
where $\Gamma'_j=E\left(\Gamma_j,v\right)$.
}

\medskip
Under the conditions of Lemma 2.3 we define $C_K(\Gamma, \Gamma_1)$ as $C_K(\Gamma, \Gamma_1):=\Gamma_2$ and $C_D(\Gamma, \Gamma_1):=E\left(\Gamma_2,v\right)$. In particular, if $\Gamma_1$ is the graph consisting of one edge $r\in\mathcal{E}\cup\mathcal{R}$ then we denote the obtaining graphs as $C_K(\Gamma, r)$ and $C_D(\Gamma, r)$ and say that $C_K(\Gamma, r)$ is the graph obtaining from $\Gamma$ by $K$-cutting-off the edge $r$ and $C_D(\Gamma, r)$ is the graph obtaining from $\Gamma$ by $D$-cutting-off $r$.

\medskip
{\bf Corollary 2.2.} {\it Let $r=[v,\infty)\in\mathcal{R}(\Gamma)$,
$v\in V(\Gamma)$. Then
$$
\Delta(\lambda)=d_r(\rho)\Delta^r(\lambda)+d^r(\rho)\Delta_r(\lambda),
$$
where $\Delta(\lambda)=\Delta\left(\lambda, \Gamma\right)$, $\Delta^r(\lambda)=\Delta\left(\lambda, C_K(\Gamma,r)\right)$,
$\Delta_r(\lambda)=\Delta\left(\lambda, C_D(\Gamma, r)\right)$, $d_r(\rho)=e_r(v,\rho)$,  $d^r(\rho)=\partial_r e_r(v,\rho)$.}

\medskip
{\bf Remark 2.1.} It is often convenient to use both spectral
parameters $\lambda$ and $\rho$ in the same formula like it has been
done in Corollary 2.2. Here and everywhere below we assume
$\lambda=\rho^2$ and if $\rho\in\mathbf{R}\setminus\{0\}$ we agree that
$\lambda=\rho^2+\mbox{sgn}\rho\cdot i0$ (i.e. on the boundary of the cut in
$\mathbf{C}\setminus [0,+\infty)$ one should take here and below the
corresponding limit).

\medskip
{\bf Corollary 2.3.} {\it Let $r\in\mathcal{E}(\Gamma)$ be the edge connecting the vertices $u$ and $v$, where $u\in\partial\Gamma$. Then
$$
\Delta(\lambda)=d_r(\lambda)\Delta^r(\lambda)+d^r(\lambda)\Delta_r(\lambda),
$$
where $\Delta(\lambda)=\Delta\left(\lambda, \Gamma\right)$, $\Delta^r(\lambda)=\Delta\left(\lambda, C_K(\Gamma,r)\right)$,
$\Delta_r(\lambda)=\Delta\left(\lambda, C_D(\Gamma, r)\right)$, $d^r(\lambda)=\Delta\left(\lambda, r^*\right)$, $d_r(\lambda)=\Delta\left(\lambda, r_*\right)$. $r^*$, $r_*$ are the graphs with one edge $r$, $V(r^*)=V(r_*)=\{u,v\}$. Further, $v$ is the $D$-type boundary vertex for $r_*$ and $K$-type boundary vertex for $r^*$. Type of boundary vertex $u$ for both $r_*$ and $r^*$ is the same as in $\Gamma$. The orderings of the vertices $u,v$ in both these graphs are the same as in $\Gamma$ as well.}

\medskip
Let us consider the zeros of the characteristic function $\Delta(\lambda,\Gamma)$. First we note that the set $\Lambda^-(\Gamma, q)$ of all negative eigenvalues of $L$ coincides with the set of all negative zeros of $\Delta(\lambda,\Gamma)$. Let $\Lambda^-_0(\Gamma, q)$ be the set of
zeros of $\Delta(\lambda, \Gamma, q)$ in $\mathbf{C}\setminus[0,+\infty)$
counted with multiplicity.
Denote
$N_{-}(\Gamma,q)=\mbox{card}(\Lambda^-_0\left(\Gamma,q)\right)$. Proceeding as in {\cite{FrI}} we arrive at the following assertion.

\medskip
{\bf Lemma 2.4.} {\it Let $\Gamma$ be an arbitrary noncompact graph. Then the following estimate holds:
$$
N_{-}(\Gamma,q)\leq N_0+Q,
$$
where
$$
Q=\sum\limits_{r\in\mathcal{R}} \int\limits_r |x|\cdot |q(x)|\cdot d|x|
$$
and $N_0$ depends only upon $q(x)$, $x\in \bigcup\limits_{r\in\mathcal{E}}r
$, i.e. upon the values of $q(\cdot)$ on the compact
part of $\Gamma$.
 }

\medskip
Now we consider the positive zeros of $\Delta(\lambda,\Gamma)$. More exactly, let $\Lambda^+_0(\Gamma)$ be the set of all positive zeros of the function $\Delta(\Gamma,\lambda+i0)$. First we need the following estimates that can be obtained in a similar way as Lemma 2.4 in {\cite{FrI}}.

{\bf Lemma 2.5.} {\it In terms of Corollary 2.2 the following estimates hold}
$$
\left|\Delta(\lambda)\right|\geq\left|\Delta_r(\lambda)\right|\cdot\left|d_r(\rho)\right|\cdot\left|
\mbox{Im}\ m_r(\lambda)\right|,
$$
$$
\left|\Delta(\lambda)\right|\geq\left|\Delta^r(\lambda)\right|\cdot\left|d^r(\rho)\right|\cdot\left|
\mbox{Im}\ \frac{1}{m_r(\lambda)}\right|,
$$
{\it where $\rho\in\overline{\Omega}_+\setminus\{0\}$,
$m_r(\lambda)=d^r(\rho)(d_r(\rho))^{-1}$ is the classical Weyl
function for $r$.}

\medskip
Now we can obtain the following result

\medskip
{\bf Lemma 2.6.} {\it $\Lambda^+_0(\Gamma)$ is at most countable set. The set $Z^+_0(\Gamma):=\{\rho: \rho^2\in\Lambda^+_0(\Gamma) \}$ has the following property: for any segment $[t,t+1]$ the number of elements of $Z^+_0(\Gamma)$ lying in this segment is bounded by some constant which does not depend on $t$.}

\medskip
{\bf Proof.} Since for any positive $\rho$ one has $d_r(\rho)\neq 0$ and $\mbox{Im}m_r(\lambda+i0)>0$ we conclude that $\Delta(\lambda+i0,\Gamma)=0$ implies $\Delta(\lambda+i0,\Gamma')=0$, where $\Gamma'=C_K(\Gamma,r)$. We can repeat this and cut-off subsequently all the rays. Thus, any $\lambda_0\in\Lambda^+_0(\Gamma)$ must be a zero of $\Delta(\lambda,\Gamma_c)$, where $\Gamma_c$ is a compact graph obtained from $\Gamma$ by cutting-off all the rays. For compact graphs validity of assertion of the Lemma is well-known {\cite{Yur5}}. $\hfil\Box$

\medskip
Now we consider the characteristic function $\Delta(\lambda, G)$ of (arbitrary) A-graph $G$. Denote $|r|$ the length of the edge $r\in\mathcal{E}(G)$ and $|G|:=\sum\limits_{r\in\mathcal{E}(G)}|r|$. Define the set $\mathcal{E}^{\pm}=\{\sum\limits_{r\in\mathcal{E}}\varepsilon_r|r|: \varepsilon_r\in\{-1,0,1\}\}$.

\medskip
{\bf Lemma 2.7.} {\it For $\lambda=\rho^2$, $\rho\to\infty$, $\rho\in\overline{\Omega}_+$ the following asymptotical representation holds:
$$
\Delta(\lambda, G)=\left(\frac{i}{2\rho}\right)^{N(G)-1}\left(\sum\limits_{l\in\mathcal{E}^{\pm}}B_l(G)\exp(-i\rho l)+O(\rho^{-1}\exp(\tau|G|))\right),
$$
where $N(G)=N_D(G)+N_{\mathcal{C}}(G)$, $N_D(G)$ is the number of $D$-type boundary vertices, $N_{\mathcal{C}}(G)$ is the number of cycles, $\tau=\mbox{Im}\rho$ and $B_l(G)$ are the constants that do not depend upon the potential $q(\cdot)$. Moreover, all the $B_l(G)$, $l\in\mathcal{E}^{\pm}$ are real and $B_{|G|}(G)\neq 0$.
 }

\medskip
{\bf Proof.} We use the induction with respect to the number of edges. For any one-edge graph (i.e. graph consisting of one simple edge or one one-edge cycle) the required assertion can be obtained via the direct calculation. Now we assume the assertion to be true for any A-graph with less than $n$ edges and consider an arbitrary A-graph $G$ with $n$ edges. Let us take some internal vertex $v$ such that $I(v)$ contains at least 3 edges (if we could not find such vertex the situation is actually equivalent to the case of one-edge graph mentioned above). Then we can represent $G$ as $G=\bigcup\limits_{k=1}^p G_k$, where:
\begin{itemize}
\item all $G_j$ are A-graphs with less then $n$ edges;
\item for any $j\neq k$ $G_j\cap G_k=v$;
\item if $v\in\partial G_j$ then $v\in\partial_K G_j$;
\item each $G_j$ has exactly 1 $a$-edge containing vertex $v$.
\end{itemize}
The last requirement guarantees, in particular, that all $G'_j:=E(G_j,v)$ are A-graphs as well. Thus we can use the representation from Corollary 2.1 for $\Delta(\lambda, G)$ and (by the inductive assumption) assertion of Lemma for each of $\Delta(\lambda, G_j)$, $\Delta(\lambda, G'_j)$.

Let us consider the values of $N(G_j)$ and $N(G'_j)$. Note that $v$ is either boundary vertex for $G_j$ or the vertex belonging to some cycle of $G_j$. In first case we have $N_\mathcal{C}(G'_j)=N_\mathcal{C}(G_j)$, $N_D(G'_j)=N_D(G_j)+1$. In second case we have $N_\mathcal{C}(G'_j)=N_\mathcal{C}(G_j)-1$, $N_D(G'_j)=N_D(G_j)+2$ and in both cases we obtain $N(G'_j)=N(G_j)+1$. Since $\sum\limits_{j=1}^p N(G_j)=N(G)$ the calculation described above yields the required representation for $\Delta(\lambda, G)$ with the constants $B_l(G)$ that are real and independent of $q(\cdot)$ (because this was true for all $B_l(G_k)$, $B_l(G'_k)$). Now we are to control the value of $B_{|G|}(G)$. Simple algebra yields:
$$
B_{|G|}(G)=\sum\limits_{k=1}^p B_{|G_k|}(G_k)\prod\limits_{j\neq k} B_{|G_j|}(G'_j).
$$
By the inductive assumption we have $B_{|G_j|}(G'_j)\neq 0$ and we can rewrite the last relation as follows:
$$
B_{|G|}(G)=\prod\limits_{j=1}^p B_{|G_j|}(G'_j)\sum\limits_{k=1}^p \frac{B_{|G_k|}(G_k)}{B_{|G_k|}(G'_k)}. \eqno(2.7)
$$
Let us consider the Weyl functions $M_v(\lambda,G_j)$. The representation from Lemma yields the following asymptotics for $\rho\to \infty$, $0<\alpha<\mbox{arg}\rho<\beta<\pi/2$:
$$
M_v(\lambda,G_j)=\frac{\Delta(\lambda, G_j)}{\Delta(\lambda, G'_j)}=-2i\rho\frac{B_{|G_j|}(G_j)}{B_{|G_j|}(G'_j)}(1+o(1)).
$$
Since $M_v(\lambda,G_j)$ are Nevanlinna functions we conclude that all $B_{|G_j|}(G_j)\left(B_{|G_j|}(G'_j)\right)^{-1}$ are real and negative. This means that the sum in right-hand side of (2.7) is nonzero and consequently $B_{|G|}(G)\neq 0$. $\hfil\Box$

\medskip
Let us agree to use the notation $A_\varepsilon$, $\varepsilon>0$ for (different) sets of the form $A_\varepsilon=\{\rho\in\overline\Omega_+:\mbox{dist}(\rho,Z)\geq \varepsilon\}$, where $Z \subset\{\rho: 0\leq\mbox{Im}\rho \leq\tau_0\}$
 is some at most countable set with the property: for any real $t$ the number of elements of $Z$ lying in the rectangle $\{\mbox{Re}\rho\in[t,t+1], \mbox{Im}\rho\in[0,\tau_0]\}$ is bounded by some constant which does not depend on $t$.

From Lemma 2.7 using standard methods {\cite{BeKu}} one can deduce the following result.

\medskip
{\bf Corollary 2.4.} {\it For $|\rho|>\rho_*$, $\rho\in A_\varepsilon$ the following estimates hold:
$$
C_1|\rho|^{1-N(G)}\exp(\tau|G|)<\left|\Delta(\lambda,G)\right|<C_2|\rho|^{1-N(G)}\exp(\tau|G|).
$$
}

\medskip
{\bf 3. Particular inverse scattering problem on the ray.} Let us take an arbitrary ray
$r\in\mathcal{R}(G)$. We call the function
$\psi_r(x,\rho)$, $x\in G$, $\rho\in
\Omega_{+}$ {\it the Weyl--type solution
associated with $r$} iff:\\[5mm]
1) it is continuous on $G$ (with respect to $x$) and satisfying $MC(v)$ for all $v\in V$; \\
2) it solves the differential equation $\ell \psi_r=\rho^2 \psi_r$, $x\in
\mbox{int}\
r'$, $r'\in \mathcal{E}(G)\cup\mathcal{R}(G)$; \\
3) $\psi_r(x,\rho)=O\left(\exp(i\rho|x|)\right)$ as $x\to\infty$, $x\in r'$,
$r'\in\mathcal{R}\setminus\{r\}$;\\
4) $\psi_r(x,\rho)=\exp(-i\rho|x|)(1+o(1))$ as $x\to\infty$, $x\in r$.\\

Proceeding in a similar way as in {\cite{FrI}} one can obtain the following results.

\medskip
{\bf Lemma 3.1.} {\it For $x\in r$ $\psi_r(x,\rho)$ is
meromorphic with respect to $\rho$ in $\Omega_+$ with possible poles
on the imaginary axis.}

\medskip
We denote the set of poles of $\psi_r(x,\rho)$,
$x\in r$ as $Z_r^-$.

\medskip {\bf Lemma 3.2.} {\it
$Z_r^-$ is a finite set. If $\rho_0\in Z_r^-$ then
$\lambda_0=\rho_0^2\in \Lambda^-$.}

\medskip {\bf Lemma 3.3.} {\it All poles of $\psi_r(x,\rho)$, $x\in r$ are simple.
For the residue ${\rm res}_{\rho=\rho_0} \psi_r(x,\rho)$,
$\rho_0\in Z_r^-$ the following representation holds:
$$
{\rm res}_{\rho=\rho_0} \psi_r(x,\rho)=i\alpha_r(\rho_0)e_r(x,\rho_0).
$$
The values $\alpha_r(\rho_0)$ are all real and positive.}

\medskip
We call the values $\alpha_r(\rho_0)$, $\rho_0\in Z^-_r$ {\it the weight
numbers}.

\medskip
Denote $Z_0^+$ the
set of all $\rho\in\mathbf{R}$ such that $\lambda=\rho^2\in\Lambda_0^+:=\Lambda_0^+(G)$.

\medskip
{\bf Lemma 3.4.} {\it If $\rho_0
\in\mathbf{R}\setminus\left(\{0\}\cup Z_0^+\right)$ then  there exists the limit $\psi_r(x,\rho_0):=\lim\limits_{\rho\to\rho_0,
\rho\in\Omega_+}\psi_r(x,\rho)$. If $\rho_0\in Z_0^+$ then
$\psi_r(x,\rho)$ and $\psi'_r(x,\rho)$ are bounded as
$\rho\to\rho_0, \rho\in\Omega_+$.}

\medskip
{\bf Lemma 3.5.} {\it For $\psi_r(x,\rho)$,
$\rho\in\mathbf{R}\setminus\left(\{0\}\cup Z_0^+\right)$ the
following representation holds:
$$
\psi_r(x,\rho)=e_r(x,-\rho)+s_r(\rho)e_r(x,\rho), \quad
x\in r.
$$}

\medskip
We call the function $s_r(\cdot)$,
{\it the
reflection coefficient associated with $r$}.

\medskip
{\bf Lemma 3.6.} {\it For all
$\rho\in\mathbf{R}\setminus\left(\{0\}\cup Z_0^+\right)$ one has $s_r(-\rho)=\overline{s_r(\rho)}$ and
$|s_r(\rho)|\leq 1$.}

\medskip
{\bf Lemma 3.7.} {\it $\psi_r(x,\rho)$, $\psi'_r(x,\rho)$,
$x\in r$ are bounded as $\rho\to 0$,
$\rho\in\overline\Omega_+$.}

\medskip
Now we agree that together with $L=L(q,G)$ we consider an operator
$\tilde L=L(G, \tilde q)$ on the same graph $G$ but having a different
potential $\tilde q(\cdot)$ satisfying the same conditions as $q(\cdot)$. If a certain symbol $\xi$ denotes an object
related to $L$, then the corresponding symbol $\tilde\xi$ with tilde
denotes the analogous object related to $\tilde L$ and
$\hat\xi:=\xi-\tilde\xi$.

\medskip
{\bf Lemma 3.8.} {\it For $x\in r$, $\rho\to\infty$,
$\rho\in A_\varepsilon$ the following estimates hold:
$$
\psi_r(x,\rho)=O\left(\exp(-i\rho|x|)\right), \
\psi'_r(x,\rho)=O\left(\rho\exp(-i\rho|x|)\right),
$$
$$
\hat\psi_r(x,\rho)=O\left(\rho^{-1}\exp(-i\rho|x|)\right).
$$
}

\medskip
{\bf Proof.} In order to obtain the asymptotics for $\psi_r(x,\rho)$ it is convenient to use the following representation, that can be obtained by direct calculation:
$$
\psi_r(x,\rho)=\gamma_r(\rho)e_r(x,\rho)+\delta_r(\rho)S_r(x,\lambda),
\quad x\in r, \eqno(3.1)
$$
where
$$
\delta_r(\rho)=-\frac{2i\rho}{d_r(\rho)}, \eqno(3.2)
$$
$$
\gamma_r(\rho)=\frac{2i\rho}{d_r(\rho)}\cdot
\frac{\Delta_r(\lambda)}{\Delta(\lambda)} \eqno(3.3)
$$
and $\Delta_r(\lambda)$ is the characteristic function for
$G_r:=C_D(G,r)$ (we recall that $d_r(\rho)=e_r(v,\rho)$).

First we estimate $\gamma_r(\rho)$. Using Corollary 2.4 and taking into account that $|G_r|=|G|$ and $N(G_r)\geq N(G)+1$ we obtain
$$
\frac{\Delta_r(\lambda)}{\Delta(\lambda)}\leq \frac{C}{|\rho|}
$$
that yields
$$\gamma_r(\rho)=O(1), \ \rho\to\infty, \ \rho\in A_\varepsilon. \eqno(3.4)$$
Now consider $\hat\gamma_r(\rho)$. From Lemma 2.7 and Corollary 2.4 one can deduce the following estimates that hold for $|\rho|>\rho_*$, $\rho\in A_\varepsilon$:
$$
\frac{\hat\Delta(\lambda)}{\Delta(\lambda)}=O\left(\frac{1}{\rho}\right), \ \frac{\hat\Delta_r(\lambda)}{\Delta_r(\lambda)}=O\left(\frac{1}{\rho}\right).
$$
This yields
$$\hat\gamma_r(\rho)=O\left(\rho^{-1}\right) \eqno(3.5)$$
for $|\rho|>\rho_*$, $\rho\in A_\varepsilon $.

To complete the proof it is sufficient to use the estimates (3.4), (3.5), the obvious estimates:
$$
\delta_r(\rho)=O(\rho), \ \hat\delta_r(\rho)=O(1)
$$
and the classical asymptotics:
$$
e^{(\nu)}_r(x,\rho)=(i\rho)^\nu\mbox{e}^{i\rho|x|}\left(1+O(\rho^{-1})\right), \
\hat e_r(x,\rho)=O\left(\rho^{-1}\mbox{e}^{i\rho|x|}\right),
$$
$$
S^{(\nu)}_r(x,\lambda)=O\left(\rho^{\nu-1}\mbox{e}^{-i\rho|x|}\right), \ \hat
S_r(x,\lambda)=O\left(\rho^{-2}\mbox{e}^{-i\rho|x|}\right).
$$
$\hfil\Box$

\medskip
{\bf Definition 3.1.} The data $J_r:=\{s_r(\cdot), Z^-_r,
\alpha_r(\rho), \rho\in Z_r^- \}$ are called {\it the scattering
data, associated with $r$}.

\medskip
{\bf Problem IP1(r).} Given $J_r$, recover the potential $q(x)$ for $x\in r$.

\medskip
{\bf Theorem 3.1.} {\it If $J_r=\tilde J_r$ then $q=\tilde q$ a.e. on $r$, i.e. the potential on the ray $r$ is uniquely determined by the  scattering data, associated with $r$. Moreover, $M_v(\cdot, G)=\tilde M_v(\cdot, G)$.}

\medskip
{\bf Proof.} Consider for $x\in r$,
$\lambda\in\mathbf{C}\setminus[0,+\infty)$ the following functions:
$$
\varphi_1(x,\lambda):=\psi_r(x,\rho), \quad
\varphi_2(x,\lambda):=e_r(x,\rho), \quad \lambda=\rho^2,
\rho\in\Omega_+.
$$

Let us define the matrices
$$
\Psi(x,\lambda):= \left[
\begin{array}{ll}
\varphi_1(x,\lambda) & \varphi_2(x,\lambda)\\
\varphi'_1(x,\lambda) & \varphi'_2(x,\lambda)
\end{array}
\right]
$$
and $\tilde\Psi(x,\lambda)$ and introduce the {\it spectral mapping
matrix}:
$$
P(x,\lambda):=\Psi(x,\lambda)\tilde\Psi^{-1}(x,\lambda).
$$

It follows from Lemma 3.5 that for the limit-value matrices
$\Psi^{\pm}(x,\lambda):=\Psi(x,\lambda\pm i0)$,
$\lambda\in(0,+\infty)\setminus\Lambda^+_0$ the following relation
holds:
$$
\Psi^-(x,\lambda)=\Psi^+(x,\lambda)w(\lambda),
$$
where
$$
w(\lambda)= \left[
\begin{array}{ll}
\overline {s_r(\rho)} & 1\\
1-|s_r(\rho)|^2 & -s_r(\rho)
\end{array}
\right],\quad \lambda=\rho^2, \rho\in(0,+\infty).
$$

Suppose that $s_r=\tilde s_r$. Then $w=\tilde w$ and
consequently $P^+(x,\lambda)=P^-(x,\lambda)$,
$\lambda\in(0,+\infty)\setminus(\Lambda^+_0\cup\tilde\Lambda^+_0)$. This
means that $P(x,\lambda)$ is holomorphic in
$\lambda\in\mathbf{C}\setminus\left(\{0\}\cup\Lambda^+_0\cup\tilde\Lambda^+_0\cup\Lambda^-_r\cup\tilde\Lambda^-_r\right)$, where $\Lambda^-_r=\{\lambda=\rho^2, \rho\in Z^-_r\}$.
Take an arbitrary
$\lambda_0\in(0,+\infty)\cap(\Lambda^+_0\cup\tilde\Lambda^+_0)$. It follows
from Lemma 3.4 that $P(x,\lambda)$ is bounded in the neighborhood of
$\lambda_0$, so $\lambda_0$ is a removable singularity for
$P(x,\lambda)$.

Then, $J_r=\tilde J_r$ means in particular
that $Z^-_r=\tilde Z^-_r$. Taking an arbitrary
$\lambda_0=\rho_0^2, \ \rho_0\in Z^-_r$ we can conclude that $\lambda_0$ is either a
 pole or a removable singularity for $P(x,\lambda)$.
Let us consider the functions $P_{11}(x,\lambda)$ and
$P_{12}(x,\lambda)$. One has:
$$
P_{11}(x,\lambda)=\frac{1}{2i\rho}\left(\psi_r(x,\rho)\tilde
e'_r(x,\rho)-\tilde\psi'_r(x,\rho)e_r(x,\rho)\right),
$$
$$
P_{12}(x,\lambda)=\frac{1}{2i\rho}\left(\tilde\psi_r(x,\rho)
e_r(x,\rho)-\psi_r(x,\rho)\tilde e_r(x,\rho)\right).
$$
Substituting here the representations
$$
\psi_r(x,\rho)=\frac{i\alpha_r(\rho_0)}{\rho-\rho_0}e_r(x,\rho_0)+O(1),
\quad\rho\to\rho_0,
$$
$$
\tilde\psi_r(x,\rho)=\frac{i\tilde\alpha_r(\rho_0)}{\rho-\rho_0}\tilde
e_r(x,\rho_0)+O(1), \quad\rho\to\rho_0,
$$
and taking into account that
$\alpha_r(\rho_0)=\tilde\alpha_r(\rho_0)$ we obtain
$P_{11}(x,\lambda)=O(1)$, $P_{12}(x,\lambda)=O(1)$ in a neighborhood
of $\lambda_0$. Thus $\lambda_0$ is a removable singularity.

Then, using Lemma 3.8 and the classical asymptotics for the Jost
solution $e_r(x,\rho)$, one can obtain the estimates:
$$
P_{11}(x,\lambda)-1=O\left(\frac{1}{\rho}\right), \quad
P_{12}(x,\lambda)=O\left(\frac{1}{\rho}\right), \quad
\lambda\to\infty, \rho^2=\lambda, \rho\in A_\varepsilon.
$$
On the other hand Lemma 3.7 yields:
$$
P_{11}(x,\lambda)-1=O\left(\frac{1}{\rho}\right), \quad
P_{12}(x,\lambda)=O\left(\frac{1}{\rho}\right), \quad \lambda\to 0,
\rho^2=\lambda.
$$
These estimates together mean that actually $P_{11}(x,\lambda)-1=0$,
$P_{12}(x,\lambda)=0$, i.e.
$\varphi_\nu(x,\lambda)=\tilde\varphi_\nu(x,\lambda)$, $\nu=1,2$
and, consequently, $q(x)=\tilde q(x)$ for a.e. $x\in r$. Notice, that
in particular we have $\psi_r(x,\rho)=\tilde\psi_r(x,\rho)$, $x\in r$, $\rho\in\overline\Omega_+\setminus Z^-_r$.

Since we have $\psi_r(x,\rho)=\psi_r(v,\rho)\cdot\Phi_v(x,\lambda,G^r)$, $x\in G^r:=C_K(G,r)$ the matching condition $MC(v)$ for $\Phi_v(x,\lambda,G^r)$ reads as follows:
$$
\frac{\partial_r\psi_r(v,\rho)}{\psi_r(v,\rho)}+M_v(\lambda,G^r)=0
$$
and we obtain $M_v(\lambda,G^r)=\tilde M_v(\lambda,G^r)$. Finally, since we have
$$
M_v(\lambda,G)=M_v(\lambda,G^r)+m_r(\lambda)
$$
we can conclude now that $M_v(\lambda,G)=\tilde M_v(\lambda,G)$.
$\hfil\Box$

\medskip
{\bf 4. Particular inverse spectral problem for compact boundary edge.} Let us consider some edge $r\in\mathcal{E}$ connecting the vertices $u$ and $v$, where $v$ is a boundary vertex.

\medskip
{\bf Problem IP2(r).} Given the Weyl function $M_v(\cdot, G)$, recover the potential $q(x)$ for $x\in r$.

In our studying this problem we follow the standard scheme of the spectral mapping method {\cite{Yur4}}, {\cite{FrYu01}}. First we need some asymptotics for Weyl solution $\Phi_v(x,\lambda)$, $x\in r$.

\medskip
{\bf Lemma 4.1.} {\it For $\lambda=\rho^2$, $\rho\to\infty$, $\rho\in A_\varepsilon$ with any $\varepsilon>0$ the following asymptotics hold:
$$
\Phi_v(x,\lambda)=O\left(\exp(-\tau|x-v|)\right), \ \Phi'_v(x,\lambda)=O\left(\rho \exp(-\tau|x-v|)\right),
$$
$$
\hat\Phi_v(x,\lambda)=O\left(\rho^{-1}\exp(-\tau|x-v|)\right),
$$
where $\tau=\mbox{Im}\rho$ and the derivative $\Phi'_v(x,\lambda)$ is considered with respect to the natural parameter measured along the edge $r$ from the vertex $v$.
}

\medskip
{\bf Proof.} For definiteness we assume that $v$ is of $D$-type (otherwise the representations below using the characteristic functions require slight modifications but the result remains the same).

We use the representation:
$$
\Phi_v(x,\lambda)=\gamma_r(\lambda)S_{r,v}(x,\lambda)+\delta_r(\lambda)S_{r,u}(x,\lambda),\  x\in r, \eqno(4.1)
$$
where $S_{r,v}(x,\lambda)$, $S_{r,u}(x,\lambda)$ are the (local) solutions for the equation $\ell y=\lambda y$ on the edge $r$ normalized by the initial conditions: $S_{r,v}(v,\lambda)=S_{r,u}(u,\lambda)=0$, $\partial_r S_{r,v}(v,\lambda)=\partial_r S_{r,u}(u,\lambda)=1$. Direct calculation yields the following representations for the coefficients $\gamma_r(\lambda)$, $\delta_r(\lambda)$:
$$
\gamma_r(\lambda)=-\frac{1}{d_r(\lambda)}\frac{\Delta_r(\lambda)}{\Delta(\lambda)}, \ \delta_r(\lambda)=\frac{1}{d_r(\lambda)}, \eqno(4.2)
$$
where we use the same notations as in Corollary 2.3.

First we estimate $\gamma_r(\lambda)$. Using Corollary 2.4 and taking into account that $|G_r|=|G|-|r|$ and $N(G_r)\geq N(G)$ we obtain
$$
\frac{\Delta_r(\lambda)}{\Delta(\lambda)}\leq C\exp(-\tau|r|).
$$
Together with the classical estimate for $d_r(\lambda)=S_{r,v}(u,\lambda)$:
$$
|d_r(\lambda)|\geq C|\rho|^{-1}\exp(\tau|r|)
$$
that yields
$$\left|\gamma_r(\lambda)\right|\leq C|\rho|\exp\left(-2\tau|r|\right). \eqno(4.3)$$

Now consider $\hat\gamma_r(\lambda)$. From Lemma 2.7 and Corollary 2.4 one can deduce the following estimates that hold for $|\rho|>\rho_*$, $\rho\in A_\varepsilon $:
$$
\frac{\hat\Delta(\lambda)}{\Delta(\lambda)}=O\left(\frac{1}{\rho}\right), \ \frac{\hat\Delta_r(\lambda)}{\Delta_r(\lambda)}=O\left(\frac{1}{\rho}\right).
$$
This yields
$$\left|\hat\gamma_r(\lambda)\right|\leq C\exp\left(-2\tau|r|\right). \eqno(4.4)$$
for $|\rho|>\rho_*$, $\rho\in A_\varepsilon$.

Next, for $\delta_r$, $\hat\delta_r$ we obtain from (4.2) and classical asymptotics the following estimates:
$$
\left|\delta_r(\lambda)\right|\leq C|\rho|\exp(-\tau|r|), \ \left|\hat\delta_r(\lambda)\right|\leq C\exp(-\tau|r|). \eqno(4.5)
$$

In order to complete the proof it is sufficient now to use the representation (4.1), estimates (4.4), (4.5) and the following classical asymptotics for the local solutions:
$$
S_{r,v}(x,\lambda)=\rho^{-1}\sin\rho|x-v|+O(\rho^{-2}\exp(\tau|x-v|)),$$$$ \ S'_{r,v}(x,\lambda)=\cos\rho|x-v|+O(\rho^{-1}\exp(\tau|x-v|)),
$$
$$
S_{r,u}(x,\lambda)=\rho^{-1}\sin\rho|x-u|+O(\rho^{-2}\exp(\tau|x-u|)), $$$$\ S'_{r,u}(x,\lambda)=-\cos\rho|x-u|+O(\rho^{-1}\exp(\tau|x-u|)).
$$
$\hfil\Box$

\medskip
{\bf Theorem 4.1.} {\it If $M_v(\cdot, G)=\tilde M_v(\cdot, G)$ then $q=\tilde q$ a.e. on $r$. Moreover, $M_u(\cdot, G)=\tilde M_u(\cdot, G)$.}

\medskip
{\bf Proof.} Proceeding as in proof of Theorem 3.1 with the conventional arguments of spectral mapping method we define the matrices:
$$
\Psi(x,\lambda):= \left[
\begin{array}{ll}
\Phi_v(x,\lambda) & S_{r,v}(x,\lambda)\\
\Phi_v'(x,\lambda) & S_{r,v}'(x,\lambda)
\end{array}
\right]
$$
and $\tilde\Psi(x,\lambda)$ and introduce the  spectral mappings
matrix:
$$
P(x,\lambda):=\Psi(x,\lambda)\tilde\Psi^{-1}(x,\lambda), x\in r.
$$
Here, as in previous Lemma the derivatives are considered with respect to the natural parameter measured along the edge $r$ from the vertex $v$.

Using the representations
$$
P_{11}(x,\lambda)=\Phi_v(x,\lambda)\tilde S_{r,v}'(x,\lambda)-\tilde\Phi'_v(x,\lambda)S_{r,v}(x,\lambda),
$$
$$
P_{12}(x,\lambda)=\tilde\Phi_v(x,\lambda) S_{r,v}(x,\lambda)-\Phi_v(x,\lambda)\tilde S_{r,v}(x,\lambda),
$$
and Lemma 4.1 we obtain the estimates:
$$
P_{11}(x,\lambda)-1=O\left(\rho^{-1}\right), \ P_{12}(x,\lambda)=O\left(\rho^{-1}\right), \ \rho\to\infty, \rho\in A_\varepsilon. \eqno(4.6)
$$
On the other hand from the same representations and $M_v(\cdot, G)=\tilde M_v(\cdot, G)$ it follows that $P_{11}(x,\lambda)-1$ and $P_{12}(x,\lambda)$ are entire functions with respect to $\lambda$. In view of (4.6) we conclude that actually $P_{11}(x,\lambda)-1\equiv 0$ and $P_{12}(x,\lambda)\equiv 0$. Thus, we have $\Phi_v(x,\lambda)=\tilde\Phi_v(x,\lambda)$ and consequently $q=\tilde q$ a.e. on $r$.

Further, it is clear that $\Phi_v(x,\lambda,G)=\Phi_v(u,\lambda,G)\cdot\Phi_u(x,\lambda,G^r)$, $x\in G^r:=C_K(G,r)$. Thus the matching condition $MC(u)$ for $\Phi_v(x,\lambda, G)$ reads as follows:
$$
\frac{\partial_r \Phi_v(u,\lambda,G)}{\Phi_v(u,\lambda,G)}+M_u(\lambda,G^r)=0
$$
and we obtain $M_u(\lambda,G^r)=\tilde M_u(\lambda,G^r)$.
On the other hand the same considerations yields
$$
M_u(\lambda,G)=M_u(\lambda, r^*)+M_u(\lambda,G^r)
$$
($r^*$ is the same one-edge graph as in Corollary 2.3). Since (as it has been already proven) $q\left|\right._r=\tilde q\left|\right._r$  we obtain finally:
$
M_u(\lambda,G)=\tilde M_u(\lambda,G).
$
$\hfil\Box$

\medskip
{\bf 5. Particular inverse spectral problem for internal simple edge.} Let $r$ be internal simple edge connecting the vertices $u$ and $v$, where $u$ is nearer to the root than $v$.

\medskip
{\bf Problem IP3(r).} Given the Weyl function $M_{v}(\cdot, G)$, and $q\left|_{G^+(r)\setminus r}\right.$, recover $q\left|_r\right.$.

\medskip
{\bf Theorem 5.1.} {\it If $M_{v}(\cdot, G)=\tilde M_{v}(\cdot, G)$ and $q\left|_{G^+(r)\setminus r}\right.=\tilde q\left|_{G^+(r)\setminus r}\right.$ then $q\left|_r\right.=\tilde q\left|_r\right.$. Moreover, $M_u(\cdot, G)=\tilde M_u(\cdot, G)$.}

\medskip
{\bf Proof.} Define $G_0^+(r):=C_K\left(G^+(r), r\right)$, $G_0^-(r):=C_K\left(G, G^+_0(r)\right)$. Since
$$
M_v(\lambda,G)=M_v(\lambda,G^+_0(r))+M_v(\lambda,G^-_0(r)),
$$
under the conditions of Theorem we have $M_v(\lambda,G^-_0(r))=\tilde M_v(\lambda,G^-_0(r))$ that by virtue of Theorem 4.1 yields $q\left|_r\right.=\tilde q\left|_r\right.$. This means, in turn that $q\left|_{G^+(r)}\right.=\tilde q\left|_{G^+(r)}\right.$ and
$$
M_u(\lambda,G^+_0(r))=\tilde M_u(\lambda,G^+_0(r)). \eqno(5.1)
$$
Further, $M_v(\lambda,G^-_0(r))=\tilde M_v(\lambda,G^-_0(r))$ implies
$$\Phi_v(x,\lambda, G^-_0(r))=\tilde \Phi_v(x,\lambda, G^-_0(r)). \ x\in r \eqno(5.2)$$
Notice that the matching condition $MC(u)$ for $\Phi_v(x,\lambda, G^-_0(r))$ reads as follows:
$$
\frac{\partial_r \Phi_v(u,\lambda, G^-_0(r))}{\Phi_v(u,\lambda, G^-_0(r))}+M_u(\lambda, G^-(r))=0.
$$
In view of (5.2) this means that
$$
M_u(\lambda, G^-(r))=\tilde M_u(\lambda, G^-(r)).
$$
From this, taking into account (5.1) and the relation
$$
M_u(\lambda, G)=M_u(\lambda, G^-(r))+M_u(\lambda, G^+(r))
$$
we obtain $M_u(\lambda, G)=\tilde M_u(\lambda, G)$ and this completes the proof. $\hfil\Box$

\medskip
{\bf 6. Particular inverse spectral problem for boundary cycle.} Now we consider some boundary cycle $\mathfrak{c}\in\mathcal{C}$.

\medskip
{\bf Problem IP4$(\mathfrak{c})$.} Given the Weyl function $M_{v_\mathfrak{c}}(\cdot, G_\mathfrak{c})$, recover the potential $q\left|_\mathfrak{c}\right.$.

\medskip
{\bf Theorem 6.1.} {\it If $M_{v_\mathfrak{c}}(\cdot, G_\mathfrak{c})=\tilde M_{v_\mathfrak{c}}(\cdot, G_\mathfrak{c})$ then $q=\tilde q$ a.e. on $\mathfrak{c}$. Moreover, $M_{u_\mathfrak{c}}(\cdot, G)=\tilde M_{u_\mathfrak{c}}(\cdot, G)$.}

\medskip
{\bf Proof.} First, we can use Theorem 4.1 and conclude that $q\left|_{r'_p}\right.=\tilde q\left|_{r'_p}\right.$ and $M_{v_{p-1}}(\lambda,G_\mathfrak{c})=\tilde M_{v_{p-1}}(\lambda,G_\mathfrak{c})$. Then, using Theorem 5.1 we obtain for $j=p-1,\ldots,1$ subsequently: $q\left|_{r_j}\right.=\tilde q\left|_{r_j}\right.$ and $M_{v_{j-1}}(\lambda,G_\mathfrak{c})=\tilde M_{v_{j-1}}(\lambda,G_\mathfrak{c})$. Finally we conclude that $q\left|_{\mathfrak{c}}\right.=\tilde q\left|_{\mathfrak{c}}\right.$ and $M_{u_\mathfrak{c}}(\cdot, G_\mathfrak{c})=\tilde M_{u_\mathfrak{c}}(\cdot, G_\mathfrak{c})$.

Define $G^-(\mathfrak{c}):=C_K(G,G^+(\mathfrak{c}))$. Since
$$
M_{u_\mathfrak{c}}(\lambda,G_\mathfrak{c})=M_{u_\mathfrak{c}}(\lambda,G^+_\mathfrak{c}(r_0))+M_{u_\mathfrak{c}}(\lambda,G^-(\mathfrak{c})).
$$
and (as it has been actually proven) $q\left|_{G^+_\mathfrak{c}(r_0)}\right.=\tilde q\left|_{G^+_\mathfrak{c}(r_0)}\right.$ we have:
$$
M_{u_\mathfrak{c}}(\lambda,G^-(\mathfrak{c}))=\tilde M_{u_\mathfrak{c}}(\lambda,G^-(\mathfrak{c})).
$$
Taking into account that
$$
M_{u_\mathfrak{c}}(\lambda,G)=M_{u_\mathfrak{c}}(\lambda,G^-(\mathfrak{c}))+M_{u_\mathfrak{c}}(\lambda,G^+(\mathfrak{c}))
$$
and $q\left|_{G^+(\mathfrak{c})}\right.=\tilde q\left|_{G^+(\mathfrak{c})}\right.$ we obtain finally $M_{u_\mathfrak{c}}(\lambda,G)=\tilde M_{u_\mathfrak{c}}(\lambda,G)$. $\hfil\Box$

\medskip
{\bf 7. Particular inverse spectral problem for internal cycle.} Consider some internal cycle $\mathfrak{c}\in\mathcal{C}$.

\medskip
{\bf Problem IP5$(\mathfrak{c})$.} Given $M_{v_\mathfrak{c}}(\cdot, G_\mathfrak{c})$ and $q(x)$, $x\in G^+(\mathfrak{c})\setminus \mathfrak{c}$ , recover the potential $q\left|_\mathfrak{c}\right.$.

\medskip
{\bf Theorem 7.1.} {\it If $M_{v_\mathfrak{c}}(\cdot, G_\mathfrak{c})=\tilde M_{v_\mathfrak{c}}(\cdot, G_\mathfrak{c})$ and $q(x)=\tilde q(x)$, $x\in G^+(\mathfrak{c})\setminus \mathfrak{c}$ then $q\left|_\mathfrak{c}\right.=\tilde q\left|_\mathfrak{c}\right.$. Moreover, $M_{u_\mathfrak{c}}(\cdot, G)=\tilde M_{u_\mathfrak{c}}(\cdot, G)$.}

\medskip
{\bf Proof.} It is sufficient to repeat the arguments from the proof of Theorem 6.1. $\hfil\Box$

\medskip
{\bf 8. Global inverse scattering problem.}

\medskip
{\bf Problem IP(G).} Given $J_r$, $r\in\mathcal{R}$, $M_v(\cdot,G)$, $v\in\partial G\setminus\{v^0\}$, $M_{v_\mathfrak{c}}(\cdot, G_\mathfrak{c})$, $\mathfrak{c}\in\mathcal{C}$, recover $q(x)$, $x\in G$.

\medskip
{\bf Theorem 8.1.} {\it Problem IP(G) has at most one solution, i.e., the specified data uniquely determine the potential $q(x)$, $x\in G$.}

\medskip
{\bf Proof.} For each fixed ray $r=[v,\infty)$, $r\in \mathcal{A}^{(\omega)}$ we apply Theorem 3.1 and get $q\left|_r\right.=\tilde q\left|_r\right.$, $M_v(\cdot, G)=\tilde M_v(\cdot,G)$.

For each fixed boundary edge $r$, $r\in \mathcal{A}^{(\omega)}$ connecting vertex $v\in\partial G$ with the vertex $u$ we apply Theorem 4.1 and get $q\left|_r\right.=\tilde q\left|_r\right.$, $M_u(\cdot, G)=\tilde M_u(\cdot,G)$.

For each fixed boundary cycle $\mathfrak{c}\in \mathcal{A}^{(\omega)}$ we apply Theorem 6.1 and get $q\left|_\mathfrak{c}\right.=\tilde q\left|_\mathfrak{c}\right.$, $M_{u_\mathfrak{c}}(\cdot, G)=\tilde M_{u_\mathfrak{c}}(\cdot,G)$.

Thus, we have proved that $q\left|_\mathfrak{a}\right.=\tilde q\left|_\mathfrak{a}\right.$ for all $a$-edges $\mathfrak{a}\in\mathcal{A^{(\omega)}}$.

Fix $\mu\in\{\omega-1,\ldots,0\}$ and suppose that we have proved that $q\left|_\mathfrak{a}\right.=\tilde q\left|_\mathfrak{a}\right.$ for all $a$-edges $\mathfrak{a}\in \mathcal{A}^{(\omega)}\cup\ldots\cup\mathcal{A}^{(\mu+1)}$. Then

1) For each fixed ray $r=[v,\infty)$, $r\in \mathcal{A}^{(\mu)}$ we apply Theorem 3.1 and get $q\left|_r\right.=\tilde q\left|_r\right.$, $M_v(\cdot, G)=\tilde M_v(\cdot,G)$.

2) For each fixed boundary edge $r$, $r\in \mathcal{A}^{(\mu)}$ connecting vertex $v\in\partial G$ with the vertex $u$ we apply Theorem 4.1 and get $q\left|_r\right.=\tilde q\left|_r\right.$, $M_u(\cdot, G)=\tilde M_u(\cdot,G)$.

3) For each fixed boundary cycle $\mathfrak{c}\in \mathcal{A}^{(\mu)}$ we apply Theorem 6.1 and get $q\left|_\mathfrak{c}\right.=\tilde q\left|_\mathfrak{c}\right.$, $M_{u_\mathfrak{c}}(\cdot, G)=\tilde M_{u_\mathfrak{c}}(\cdot,G)$.

4) For each fixed internal simple edge $r$, $r\in \mathcal{A}^{(\mu)}$ connecting vertex $v\in G^+(r)$ with the vertex $u\in\partial G^+(r)$ we apply Theorem 5.1 and get $q\left|_r\right.=\tilde q\left|_r\right.$, $M_u(\cdot, G)=\tilde M_u(\cdot,G)$.

5) For each fixed internal cycle $\mathfrak{c}\in \mathcal{A}^{(\mu)}$ we apply Theorem 7.1 and get $q\left|_\mathfrak{c}\right.=\tilde q\left|_\mathfrak{c}\right.$, $M_{u_\mathfrak{c}}(\cdot, G)=\tilde M_{u_\mathfrak{c}}(\cdot,G)$.

Thus, we have proved that $q\left|_\mathfrak{a}\right.=\tilde q\left|_\mathfrak{a}\right.$ for all $a$-edges $\mathfrak{a}\in\mathcal{A^{(\mu)}}$.

Using the above--mentioned arguments successively for $\mu=\omega-1,\ldots,1,0$ we get $q=\tilde q$ a.e. on $G$. $\hfil\Box$

\medskip
{\bf Acknowledgement.} This research was supported by the Russian Fund
of Basic Research and the Taiwan National Science Council (projects
10-01-00099 and 10-01-92001-NSC-a).

\noindent Ignatyev, Mikhail\\
Department of Mathematics, Saratov University, \\
Astrakhanskaya 83, Saratov 410012, Russia, \\
e-mail: mikkieram@gmail.com

\end{document}